\input amstex
\input amsppt.sty
\magnification=\magstep1
 \hsize=37truecc
 \vsize=24truecm
\baselineskip=14truept
 \NoBlackBoxes
\def\q{\quad}
\def\qq{\qquad}
\def\mod#1{\ (\text{\rm mod}\ #1)}

\def\t{\text}
\def\qtq#1{\q\t{#1}\q}
\def\mod#1{\ (\text{\rm mod}\ #1)}
\def\qtq#1{\q\t{#1}\q}
\def\f{\frac}
\def\e{\equiv}
\def\b{\binom}

\def\sls#1#2{(\f{#1}{#2})}
 \def\ls#1#2{\big(\f{#1}{#2}\big)}
\def\Ls#1#2{\Big(\f{#1}{#2}\Big)}
\let \pro=\proclaim
\let \endpro=\endproclaim
\topmatter
\title Congruences concerning Legendre polynomials II
\endtitle
\author ZHI-Hong Sun\endauthor
\affil School of Mathematical Sciences, Huaiyin Normal University,
\\ Huaian, Jiangsu 223001, PR China
\\ E-mail: zhihongsun$\@$yahoo.com
\\ Homepage: http://www.hytc.edu.cn/xsjl/szh
\endaffil

 \nologo \NoRunningHeads

\abstract{Let $p>3$ be a prime, and let $m$ be an integer with
$p\nmid m$. In the paper we solve some conjectures of Z.W. Sun
concerning $\sum_{k=0}^{p-1}\binom{2k}k^3/m^k\pmod{p^2}$,
$\sum_{k=0}^{p-1}\binom{2k}k\b{4k}{2k}/m^k\pmod p$ and
$\sum_{k=0}^{p-1}\binom{2k}k^2\b{4k}{2k}/m^k\pmod {p^2}.$
 In particular, we show that
$\sum_{k=0}^{\frac{p-1}2}\binom{2k}k^3\e 0\pmod {p^2}$ for $p\e
3,5,6\pmod 7$.  Let  $P_n(x)$ be the Legendre polynomials. In the
paper we also show that $ P_{[\frac p4]}(t)\e
-\big(\frac{-6}p\big)\sum_{x=0}^{p-1} \big(\frac{x^3-\f
32(3t+5)x-9t-7}p\big)\pmod
 p$ and
determine $P_{\frac{p-1}2}(\sqrt 2), P_{\frac{p-1}2}(\frac{3\sqrt
2}4), P_{\frac{p-1}2}(\sqrt {-3}),P_{\frac{p-1}2}(\frac{\sqrt 3}2),
P_{\frac{p-1}2}(\sqrt {-63}), P_{\frac{p-1}2}(\frac {3\sqrt 7}8)
\pmod p$, where $t$ is a rational $p-$integer, $[x]$ is the greatest
integer not exceeding $x$ and $(\frac ap)$ is the Legendre symbol.
  As consequences we determine
 $P_{[\frac p4]}(t)\pmod p$ in the cases $t=-\frac 53,-\frac
 79,-\frac {65}{63}$ and confirm many
 conjectures of Z.W. Sun.
\par\q
\newline MSC: Primary 11A07, Secondary 33C45, 11E25, 11L10, 05A10,
05A19
 \newline Keywords:
Legendre polynomial; congruence; character sum; binary quadratic
form}
 \endabstract
  \footnote"" {The author is
supported by the National Natural Sciences Foundation of China
(grant no. 10971078).}
\endtopmatter
\document
\subheading{1. Introduction}
\par  Let $p$ be an odd prime.
Following [A] we define
$$A(p,\lambda)=\sum_{k=0}^{\f{p-1}2}\b{\f{p-1}2}k^2\b{\f{p-1}2+k}{2k}\lambda^{kp}.$$
It is easily seen ([S3, Lemmas 2.2 and 2.4]) that
$$\b{\f{p-1}2}k\e \f{\b{2k}k}{(-4)^k}\mod p\qtq{and}
\b{\f{p-1}2+k}{2k}\e \f{\b{2k}k}{(-16)^k}\mod{p^2}.\tag 1.1$$ Thus,
by Fermat's little theorem, for any rational $p$-integer $\lambda$,
$$A(p,\lambda)\e
\sum_{k=0}^{\f{p-1}2}\b{2k}k^3\Ls{\lambda}{64}^k\mod p.\tag 1.2$$
 \par For
positive integers $a,b$ and $n$, if $n=ax^2+by^2$ for some integers
$x$ and $y$, we briefly say that $n=ax^2+by^2$. In 1985, Beukers[B]
conjectured a congruence for $A(p,1)\mod{p^2}$ equivalent to
$$\sum_{k=0}^{(p-1)/2}\f{\b{2k}k^3}{64^k}\e
\cases 0\mod{p^2}&\t{if $p\e 3\mod 4$,}\\4x^2-2p\mod{p^2}&\t{if
$p=x^2+4y^2\e 1\mod 4$.}\endcases\tag 1.3$$ This congruence was
proved by several authors including Ishikawa[Is]($p\e 1\mod 4$), van
Hamme[vH]($p\e 3\mod 4$) and Ahlgren[A]. In 1998, by using the
hypergeometric series $_3F_2(\lambda)_p$ over the finite field
$F_p$, Ono[O] obtained congruences for $A(p,\lambda)\mod p$ in the
cases $\lambda=-1,-8,-\f 18,4,\f 14,64,\f 1{64}$, see also [A].
Hence from (1.2) we deduce congruences for $\sum_{k=0}^{\f{p-1}2}\f
1{m^k}\b{2k}k^3\mod p$ in the cases $m=1,-8,-64,256,-512,4096$.
Recently the author's brother Zhi-Wei Sun[Su1,Su3,Su4] conjectured
corresponding congruences modulo $p^2$. In particular, he
conjectured that
$$\aligned\sum_{k=0}^{p-1}\b{2k}k^3\e \cases 0\mod{p^2}&\t{if $p\e 3,5,6\mod 7$,}
\\4x^2-2p\mod {p^2}&\t{if $p=x^2+7y^2\e 1,2,4\mod 7$.}\endcases\endaligned\tag 1.4$$
 We note that $p\mid \b{2k}k$ for
$\f{p+1}2\le k\le p-1$. In the paper, by using Legendre polynomials
we partially solve Zhi-Wei Sun's such conjectures. For example, we
prove (1.4) in the case $p\e 3,5,6\mod 7$.

\par Let $\{P_n(x)\}$ be the Legendre polynomials given by
$$P_0(x)=1,\ P_1(x)=x,\ (n+1)P_{n+1}(x)=(2n+1)xP_n(x)-nP_{n-1}(x)\ (n\ge
1).$$
It is well known that (see [MOS, pp.\;228-232], [G,
(3.132)-(3.133)])
$$P_n(x)=\f
1{2^n}\sum_{k=0}^{[n/2]}\b nk(-1)^k\b{2n-2k}nx^{n-2k} =\f 1{2^n\cdot
n!}\cdot\f{d^n}{dx^n}(x^2-1)^n,\tag 1.5$$ where $[x]$ is the
greatest integer not exceeding $x$. From (1.5) we see that
$$P_n(-x)=(-1)^nP_n(x), \q P_{2m+1}(0)=0\qtq{and}P_{2m}(0)=\f{(-1)^m}{2^{2m}}\b{2m}m.
\tag 1.6$$ In the paper we deduce our main results for congruences
from the identity $\sum_{k=0}^n\b{2k}k^2\b{n+k}{2k}x^k$
$=P_n(\sqrt{1+4x})^2$.
\par Let $\Bbb Z$ be the set of
integers. For a prime $p$ let $\Bbb Z_p$ be the set of rational
numbers whose denominator is coprime to $p$, and let $\sls ap$ be
the Legendre symbol. On the basis of the work of Ono[O] (see also
[A, Theorem 2] and [LR]), in [S3, Theorem 2.11] the author showed
that for a prime $p>3$ and $t\in\Bbb Z_p$,
$$P_{\f{p-1}2}(t)\e
-\Ls{-6}p\sum_{x=0}^{p-1}\Ls{x^3-3(t^2+3)x+2t(t^2-9)}p\mod p.\tag
1.7$$ In the paper, using (1.7) we prove that
$$ P_{[\f p4]}(t) \e
-\Ls{6}p\sum_{x=0}^{p-1}\Ls{x^3-\f 32(3t+5)x+9t+7}p\mod
 p.\tag 1.8$$  As consequences
of (1.8), we determine $P_{[\f p4]}(-\f 53),\ P_{[\f p4]}(-\f 79),\
P_{[\f p4]}(-\f {65}{63})\mod p$ and use them to solve Z.W. Sun's
conjectures on $\sum_{k=0}^{p-1}\b{2k}k\b{4k}{2k}/m^k\mod p$ and
$\sum_{k=0}^{p-1}\b{2k}k^2\b{4k}{2k}/m^k$ $\mod {p^2}$, see [Su1].
For instance, for any prime $p>7$,
$$\sum_{k=0}^{p-1}\f{\b{2k}k^2\b{4k}{2k}}{81^k}
\e\cases 4C^2\mod p &\t{if $p\e 1,2,4\mod 7$ and so $p=C^2+7D^2$,}
\\0\mod {p^2}&\t{if $p\e 3,5,6\mod 7$.}
\endcases\tag 1.9$$

\subheading{2. Congruences for $P_{\f{p-1}2}(\sqrt x)$ and $P_{[\f
p4]}(t)\mod p$}

\pro{Lemma 2.1} Let $p$ be an odd prime and $k\in\{0,1,\ldots, [\f
p4]\}$. Then
$$\b{[\f p4]+k}{2k}\e \f 1{(-64)^k}\b{4k}{2k}\mod p\ \t{and}\  \b{p-1-2k}{\f{p-1}2}
\e \f {(-1)^{\f{p-1}2}}{16^k}\b{4k}{2k}\mod p.$$
\endpro
Proof. Suppose $r=1$ or $3$ according as $4\mid p-1$ or $4\mid p-3$.
Then clearly
$$\aligned \b{[\f p4]+k}{2k}&
=\f{(\f{p-r}4+k)(\f{p-r}4+k-1)\cdots(\f{p-r}4-k+1)}{(2k)!}
\\&\e (-1)^k\f{(4k-r)(4k-r-4)\cdots(4-r)\cdot r(r+4)\cdots(4k+r-4)}
{4^{2k}\cdot (2k)!}
\\&=\f{(-1)^k\cdot (4k)!}{2^{2k}\cdot (2k)!\cdot 4^{2k}\cdot (2k)!}
=\f{\b{4k}{2k}}{(-64)^k} \mod p
\endaligned$$
and
$$\align
&(-1)^{\f{p-1}2}\b{p-1-2k}{\f{p-1}2}\\&=(-1)^{\f{p-1}2}\f{(p-1-2k)(p-2-2k)\cdots
(p-(\f{p-1}2+2k))}{\f{p-1}2!}\\&\e \f{(2k+1)(2k+2)\cdots
(\f{p-1}2+2k)}{\f{p-1}2!} =\f{(\f{p-1}2+2k)(\f{p-1}2+2k-1)\cdots
(\f{p-1}2+1)}{(2k)!}
\\&\e \f{(4k-1)(4k-3)\cdots 3\cdot 1}{2^{2k}\cdot (2k)!}
=\f{(4k)!}{2^{2k}\cdot  (2k)!\cdot 2^{2k}\cdot (2k)!}=\f
1{2^{4k}}\b{4k}{2k}\mod p.\endalign$$ This proves the lemma.

\pro{Lemma 2.2} Let $p$ be an odd prime and let $t$ be a variable.
Then
 $$P_{[\f p4]}(t)\e \sum_{k=0}^{[p/4]}\b{4k}{2k}\b{2k}k\Ls{1-t}{128}^k\mod
 p.$$
 \endpro
 Proof. It is known that ([G, (3.135)])
 $$P_n(t)=\sum_{k=0}^n\b nk\b {n+k}k\Ls{t-1}2^k.\tag 2.1$$
  Observe that $\b nk\b{n+k}k=\b{n+k}{2k}\b
{2k}k$. By (2.1) and Lemma 2.1 we have
$$\aligned P_{[\f p4]}(t)&=\sum_{k=0}^{[p/4]}\b{[\f
p4]+k}{2k}\b{2k}k\Ls{t-1}2^k  \e
\sum_{k=0}^{[p/4]}\f{\b{4k}{2k}\b{2k}k}{(-64)^k}\Ls{t-1}2^k\mod
p.\endaligned$$ This yields the result.

 \pro{Lemma 2.3} Let $p$ be an odd prime
and let $t$ be a variable. Then
 $$ \sum_{k=0}^{[p/4]}\b{4k}{2k}\b{2k}k\Ls{1-t}{128}^k\e
(-1)^{[\f p4]}\sum_{k=0}^{[p/4]}\b{4k}{2k}\b{2k}k\Ls{1+t}{128}^k
 \mod  p.$$
 \endpro
 Proof. This is immediate from Lemma
 2.2 and (1.6).
\par\q
 \pro{Lemma 2.4} Let $p$ be an
odd prime and let $x$ be a variable. Then
$$(\sqrt
x)^{\f{p-1}2}P_{\f{p-1}2}(\sqrt x)\e (-1)^{[\f
p4]}\sum_{k=0}^{[p/4]}\f{\b{4k}{2k}\b{2k}k}{64^k}x^{\f{p-1}2-k}\mod
p.$$
\endpro
Proof. From (1.5) we see that
$$(\sqrt
x)^{\f{p-1}2}P_{\f{p-1}2}(\sqrt x)=\f
1{2^{\f{p-1}2}}\sum_{k=0}^{[p/4]}\b{\f{p-1}2}k\b{p-1-2k}{\f{p-1}2}(-1)^kx^{\f{p-1}2-k}.
$$ Thus applying (1.1) and Lemma 2.1 we obtain
$$ (\sqrt
x)^{\f{p-1}2}P_{\f{p-1}2}(\sqrt x)\e \Ls {-2}p\sum_{k=0}^{[p/4]}\f
1{(-4)^k}\b{2k}k\cdot \f 1{2^{4k}}\b{4k}{2k}(-1)^kx^{\f{p-1}2-k}\mod
p.$$ Noting that $\ls{-2}p=(-1)^{[\f p4]}$ we then obtain the
result.
 \pro{Lemma 2.5} Let $p$ be an odd prime and let $x\in \Bbb Z_p$
with $x\not\e 0,1\mod p$. Then
$$\align &\Ls{2(1-x)}p\Big(\sqrt {\f x{x-1}}\Big)^{\f{p-1}2}P_{\f{p-1}2}\Big(\sqrt {\f x{x-1}}\Big)
\\&\e (\sqrt x)^{\f{p-1}2}P_{\f{p-1}2}(\sqrt x)\e
\Ls xpP_{[\f p4]}\Big(\f 2x-1\Big)\mod p.\endalign$$
\endpro
Proof. By (1.6), Lemmas 2.2 and 2.4 we have
$$\align P_{[\f p4]}\Big(\f 2x-1\Big)&=(-1)^{[\f p4]}P_{[\f p4]}\Big(1-\f 2x\Big)
\e (-1)^{[\f p4]}\sum_{k=0}^{[p/4]}\b{4k}{2k}\b{2k}k\f 1{(64x)^k}
\\&\e \Ls xp(\sqrt x)^{\f{p-1}2}P_{\f{p-1}2}(\sqrt x)\mod
p.\endalign$$
 Therefore,
 $$\align &\Big(\sqrt {\f
x{x-1}}\Big)^{\f{p-1}2}P_{\f{p-1}2}\Big(\sqrt {\f x{x-1}}\Big)
\\&\e \Ls {x/(x-1)}pP_{[\f p4]}\Big(\f 2{x/(x-1)}-1\Big)
 =\Ls {x(x-1)}pP_{[\f p4]}\Big(1-\f 2x\Big)
\\& \e\Ls{-2(x-1)}p(\sqrt x)^{\f{p-1}2}P_{\f{p-1}2}(\sqrt x)
  \mod p.\endalign$$ This completes the proof.
 \pro{Lemma
2.6} Let $p>3$ be a prime and $u\in\Bbb Z_p$. Then
$$(\sqrt u)^{\f{p-1}2}P_{\f{p-1}2}(\sqrt u)\e -\Ls {-6}p\sum_{x=0}^{p-1}\Ls{ux^3-3(u+3)x+2(u-9)}p\mod
p.$$
\endpro
Proof. For $t\in\Bbb Z_p$ with $t\not\e 0\mod p$, by (1.7) we have
$$\align P_{\f{p-1}2}(t)&\e -\Ls {-6}p\sum_{x=0}^{p-1}\Ls{x^3-3(t^2+3)x+2t(t^2-9)}p
\\&=-\Ls {-6}p\sum_{x=0}^{p-1}\Ls{(tx)^3-3(t^2+3)tx+2t(t^2-9)}p
\\&=-\Ls {-6t}p\sum_{x=0}^{p-1}\Ls{t^2x^3-3(t^2+3)x+2(t^2-9)}p\mod
p.\endalign$$ Hence, if $u=t^2$ for some $t\in\Bbb Z_p$ with
$t\not\e 0\mod p$, then
$$\align(2\sqrt u)^{\f{p-1}2}P_{\f{p-1}2}(\sqrt
u)&=(2t)^{\f{p-1}2}P_{\f{p-1}2}(t)\e -\Ls{-3}p\sum_{x=0}^{p-1}
\Ls{ux^3-3(u+3)x+2(u-9)}p
\\&\e -\Ls{-3}p\sum_{x=0}^{p-1}((x^3-3x+2)u-9x-18)^{\f{p-1}2}\mod p.\endalign$$
This is also true for $u=0$ since $\sum_{x=0}^{p-1}\sls{-9x-18}p=0$.
Set
$$f(u)=(2\sqrt u)^{\f{p-1}2}P_{\f{p-1}2}(\sqrt
u)+\Ls{-3}p\sum_{x=0}^{p-1}((x^3-3x+2)u-9x-18)^{\f{p-1}2}.$$ Then
$f(u)\e 0\mod p$ for $u=0^2,1^2,2^2,\ldots,\sls{p-1}2^2$. From the
proof of Lemma 2.4 we know that $(2\sqrt u)^{\f{p-1}2}$
$P_{\f{p-1}2}(\sqrt u)$ is a polynomial of $u$ with degree
$\f{p-1}2$ and integral coefficients. Thus $f(u)$ is a polynomial of
$u$ with degree at most $\f{p-1}2$ and integral coefficients. As
$f(u)\e 0\mod p$ for $u=0^2,1^2,2^2,\ldots,\sls{p-1}2^2$, using
Lagrange's theorem we see that all the coefficients in $f(u)$ are
divisible by $p$. Therefore, $f(u)\e 0\mod p$ for every $u\in\Bbb
Z_p$. This yields the result.
\par Now we are ready to prove the following main result.
 \pro{Theorem 2.1} Let $p$ be
a prime greater than $3$ and $t\in\Bbb Z_p$.
\par $(\t{\rm i})$ If $t\not\e 0\mod p$, then
$$ \align (\sqrt t)^{-\f{p-1}2}P_{\f{p-1}2}(\sqrt t)
&\e P_{[\f p4]}\Ls {2-t}t\e
\sum_{k=0}^{[p/4]}\b{4k}{2k}\b{2k}k\Ls{t-1}{64t}^k \\&\e (-1)^{[\f
p4]}\sum_{k=0}^{[p/4]}\b{4k}{2k}\b{2k}k\f 1{(64t)^k}
\\& \e -\Ls{-6t}p\sum_{x=0}^{p-1}\Ls{tx^3-3(t+3)x+2(t-9)}p\mod
 p.\endalign$$
 \par $(\t{\rm ii})$ We have
  $$ P_{[\f p4]}(t)\e
\sum_{k=0}^{[p/4]}\b{4k}{2k}\b{2k}k\Ls{1-t}{128}^k \e
-\Ls{6}p\sum_{x=0}^{p-1}\Ls{x^3-\f{3(3t+5)}2x+9t+7}p\mod
 p.$$
 \endpro
Proof. By Lemmas 2.5 and 2.6 we have
$$\align&P_{[\f p4]}(2/t-1)\\&\e \Ls tp(\sqrt t)^{\f{p-1}2}
P_{\f{p-1}2}(\sqrt t) \e
-\Ls{-6t}p\sum_{x=0}^{p-1}\Ls{tx^3-3(t+3)x+2(t-9)}p\mod
p.\endalign$$
 This together with Lemmas 2.2 and 2.3 gives the
first part. Substituting $t$ by $\f 2{t+1}$ and noting that
$$\align \sum_{x=0}^{p-1}\Ls{x^3-\f{3(3t+5)}2x-9t-7}p
&=\sum_{x=0}^{p-1}\Ls{(-x)^3-\f{3(3t+5)}2(-x)-9t-7}p\\& =\Ls{-1}p
\sum_{x=0}^{p-1}\Ls{x^3-\f{3(3t+5)}2x+9t+7}p\endalign$$ we then
obtain the second part in the case $t\not\e -1\mod p$. When $t\e
-1\mod p$ we have $P_{[\f p4]}(-1)=(-1)^{[\f p4]}=\sls{-2}p$ and
$$\sum_{x=0}^{p-1}\Ls{x^3-3x-2}p=\sum_{x=0}^{p-1}\Ls{(x+1)^2(x-2)}p=\sum_{x=0}^{p-1}\Ls{x-2}p-\Ls{-1-2}p
=-\Ls{-3}p.$$ Thus the second part is also true for $t\e -1\mod p$.
The proof is now complete.
 \pro{Corollary 2.1} Let $p\ge 17$ be a prime and
$t\in\Bbb Z_p$. Then
$$\sum_{x=0}^{p-1}\Ls{x^3-\f{3(3t+5)}2x+9t+7}p=\Ls 2p
\sum_{x=0}^{p-1}\Ls{x^3+\f{3(3t-5)}2x +9t-7}p.$$
\endpro Proof. Since $P_{[\f p4]}(t)=(-1)^{[\f p4]}P_{[\f
p4]}(-t)$, by Theorem 2.3(ii) we obtain
$$\align &\sum_{x=0}^{p-1}\Ls{x^3-\f{3(3t+5)}2x+9t+7}p\\&\e
(-1)^{[\f p4]}\sum_{x=0}^{p-1}\Ls{x^3-\f{3(-3t+5)}2x -9t+7}p =\Ls
{-2}p\sum_{x=0}^{p-1}\Ls{(-x)^3+\f{3(3t-5)}2(-x)-9t+7}p \\&=\Ls
2p\sum_{x=0}^{p-1}\Ls{x^3+\f{3(3t-5)}2x+9t-7}p\mod p.\endalign$$ By
Weil's estimate ([BEW, p.183]) we have
$$\Big|\sum_{x=0}^{p-1}\Ls{x^3-\f{3(3t+5)}2x+9t+7}p\Big|\le 2\sqrt
p\qtq{and}\Big|\sum_{x=0}^{p-1}\Ls{x^3+\f{3(3t-5)}2x +9t-7}p\Big|\le
2\sqrt p.$$ Since $4\sqrt p<p$ for $p\ge 17$, from the above we
deduce the result.

\par For any prime $p>3$, in [Su1] Zhi-Wei Sun conjectured
congruences for $\sum_{k=0}^{p-1}\f{\b{4k}{2k}\b{2k}k}{m^k}$ $\mod
{p^2}$ in the cases $m=48,63,72$. Now we confirm his congruences
modulo $p$.

 \pro{Theorem 2.2} Let $p$ be a prime
greater than $3$. Then
$$\aligned P_{[\f p4]}\Big(-\f 79\Big)&\e\sum_{k=0}^{[p/4]}\f{\b{4k}{2k}\b{2k}k}{72^k}\e
(-1)^{[\f p4]}\sum_{k=0}^{[p/4]}\f{\b{4k}{2k}\b{2k}k}{576^k}
\\&\e \cases (-1)^{\f{p-1}4}\ls p32a\mod p&\t{if
$4\mid p-1$, $p=a^2+b^2$ and $4\mid a-1$,} \\0\mod p&\t{if $p\e
3\mod 4$.}\endcases\endaligned$$
\endpro
Proof.  From [BEW, Theorem 6.2.9] or [S1, (2.15)-(2.16)] we have
$$\sum_{x=0}^{p-1}\Ls {x^3-4x}p=\cases -2a&\t{if $p\e 1\mod 4$, $p=a^2+4b^2$ and $a\e 1\mod 4$,}
\\0&\t{if $p\e 3\mod 4$.}\endcases\tag 2.2$$ Thus
taking $t=9$ in Theorem 2.1(i) we deduce the result.

 \pro{Theorem 2.3} Let $p$ be a prime greater than $3$.
 Then
$$\aligned P_{[\f p4]}\Big(-\f 53\Big)&\e\sum_{k=0}^{[p/4]}\f{\b{4k}{2k}\b{2k}k}{48^k}
\e (-1)^{[\f p4]}\sum_{k=0}^{[p/4]}\f{\b{4k}{2k}\b{2k}k}{(-192)^k}
\\&\e \cases 2A\mod p&\t{if $3\mid p-1$, $p=A^2+3B^2$
 and $3\mid A-1$,}\\0\mod p&\t{if $p\e 2\mod 3$.}\endcases\endaligned$$
\endpro
Proof. It is known that (see for example [S1, (2.7)-(2.9)] or [BEW,
pp. 195-196])
$$\sum_{x=0}^{p-1}
\Ls{x^3+8}p=\cases -2A\sls 2p&\t{if $3\mid p-1$, $p=A^2+3B^2$
 and $3\mid A-1$,}\\0&\t{if $p\e 2\mod 3$.}\endcases\tag 2.3$$
 Thus, putting $t=-3$ in Theorem 2.1(i) we deduce the
 result.

\pro{Theorem 2.4} Let $p\not=2,3,7$ be a prime. Then
$$\aligned P_{[\f p4]}\Big(-\f{65}{63}\Big)&\e\sum_{k=0}^{[p/4]}\f{\b{4k}{2k}\b{2k}k}{63^k}\e
(-1)^{[\f p4]}\sum_{k=0}^{[p/4]}\f{\b{4k}{2k}\b{2k}k}{(-4032)^k}
\\&\e \cases 2C\sls p3\sls C7\mod p&\t{if $p\e 1,2,4\mod 7$ and so
 $p=C^2+7D^2$,}\\0\mod p&\t{if $p\e 3,5,6\mod
7$.}\endcases\endaligned$$
\endpro
Proof.  Putting $t=-63$ in Theorem 2.1(i) we have
$$\align &P_{[\f p4]}\Big(-\f{65}{63}\Big)\\&\e\sum_{k=0}^{[p/4]}\f{\b{4k}{2k}\b{2k}k}{63^k}\e
(-1)^{[\f p4]}\sum_{k=0}^{[p/4]}\f{\b{4k}{2k}\b{2k}k}{(-4032)^k}
\\&\e -\Ls{-6\cdot (-63)}p\sum_{x=0}^{p-1} \Ls
{-63x^3+180x-144}p=-\Ls{-42}p\sum_{x=0}^{p-1} \Ls {7x^3-20x+16}p
\\&=-\Ls{-42}p\sum_{x=0}^{p-1} \Ls {7(2x)^3-20\cdot 2x+16}p
=- \Ls{-21}p\sum_{x=0}^{p-1} \Ls {7x^3-5x+2}p\mod p.\endalign$$ From
[R1,R2] we have
$$\sum_{x=0}^{p-1}\Ls{x^3+21x^2+112x}p=\cases -2C\sls C7&\t{if $p\e 1,2,4\mod 7$ and so
 $p=C^2+7D^2$,}\\0&\t{if $p\e 3,5,6\mod
7$.}\endcases$$ As $x^3+21x^2+112x=(x+7)^3-35(x+7)-98$ we have
$$\align \sum_{x=0}^{p-1}\Ls{x^3+21x^2+112x}p&=
\sum_{x=0}^{p-1}\Ls{x^3-35x-98}p
=\sum_{x=0}^{p-1}\Ls{(-7x)^3-35(-7x)-98}p
\\&=\Ls{-1}p\sum_{x=0}^{p-1}\Ls{7x^3-5x+2}p\endalign$$
and so
$$\sum_{x=0}^{p-1}\Ls{7x^3-5x+2}p=\cases(-1)^{\f{p+1}2}2C\sls C7&\t{if $p\e 1,2,4\mod 7$ and
 $p=C^2+7D^2$,}\\0&\t{if $p\e 3,5,6\mod
7$.}\endcases\tag 2.4$$ Now combining all the above we deduce the
result.
 \pro{Theorem 2.5} Let $p\not=2,3,7$ be a prime. Then
$$\aligned&P_{\f{p-1}2}(3\sqrt 7/8)\\&\e\cases 2C\mod p&\t{if $p\e
1,9,25\mod{28}$, $p=C^2+7D^2$ and $4\mid C-1$,}
\\-2\sqrt 7D\mod p&\t{if $p\e
11,15,23\mod{28}$, $p=C^2+7D^2$ and $4\mid D-1$,}
\\0\mod p&\t{if $p\e 3,5,6\mod 7$}
\endcases\endaligned$$
and
$$\aligned&(-1)^{[\f p4]}P_{\f{p-1}2}(\sqrt {-63})\\&\e\cases 2C\mod p&\t{if $p\e
1,9,25\mod{28}$, $p=C^2+7D^2$ and $4\mid C-1$,}
\\2D\sqrt {-7}\mod p&\t{if $p\e
11,15,23\mod{28}$, $p=C^2+7D^2$ and $4\mid D-1$,}
\\0\mod p&\t{if $p\e 3,5,6\mod 7$.}
\endcases\endaligned$$
\endpro
Proof. From Lemma 2.5 and Theorem 2.4 we have
$$\aligned&\Ls{2\cdot 64}p\Ls{3\sqrt 7}8^{\f{p-1}2}P_{\f{p-1}2}\Ls{3\sqrt 7}8
\\&\e\Ls{-63}pP_{[\f p4]}\Big(-\f{65}{63}\Big)\\& \e\cases 2C\sls p3\sls
C7\mod p&\t{if $p\e 1,2,4\mod 7$ and so $p=C^2+7D^2$,}
\\0\mod p&\t{if $p\e 3,5,6\mod 7$.}
\endcases\endaligned$$
Now suppose $p\e 1,2,4\mod 7$ and so $p=C^2+7D^2$. By [S2, p.1317]
we have
$$7^{[\f p4]}\e \cases\sls C7\mod p&\t{if $p\e 1,9,25\mod{28}$ and $C\e
1\mod 4$,}
\\-\sls C7\f DC\mod p&\t{if $p\e 1,9,25\mod{28}$ and $D\e
1\mod 4$.}\endcases\tag 2.5$$
 Thus, noting that
$3^{\f{p-1}2}\ls p3\e \ls 3p\ls p3=(-1)^{\f{p-1}2}\mod p$ we get
$$\aligned &P_{\f{p-1}2}(3\sqrt 7/8)\\&\e (3\sqrt 7)^{-\f{p-1}2}2C\Ls p3\Ls C7
\\&\e\cases 7^{-\f{p-1}4}2C\Ls C7\e 2C\mod p&\t{if $p\e
1\mod 4$ and $4\mid C-1$,}
\\-\f{\sqrt 7}{7\cdot 7^{\f{p-3}4}}2C\ls C7\e \f{2\sqrt 7C^2}{7D}
\e -2\sqrt 7D\mod p&\t{if $p\e 3\mod {4}$ and $4\mid D-1$.}
\endcases\endaligned$$
Taking $x=-63$ in Lemma 2.5 we obtain $(\sqrt{-7})^{\f{p-1}2}
P_{\f{p-1}2}(\sqrt{-63})\e (\sqrt 7)^{\f{p-1}2}$
$P_{\f{p-1}2}\sls{3\sqrt 7}8$ $\mod p$. Now combining all the above
we obtain the result.

\pro{Theorem 2.6} Let $p$ be an odd prime.  Then
$$P_{\f{p-1}2}\Ls{3\sqrt 2}4\e
\cases (-1)^{\f b2}2a\mod p&\t{if $8\mid p-1$, $p=a^2+4b^2$ and
$4\mid a-1$,}
\\4b\mod p&\t{if $8\mid p-5$, $p=a^2+4b^2$
and $4\mid b-1$,}
\\0\mod p&\t{if $p\e 3\mod 4$.}
\endcases$$
\endpro
Proof. Taking $x=9$ in Lemma 2.5 we have $\ls{2\cdot(-8)}p\ls{3\sqrt
2}4^{\f{p-1}2}P_{\f{p-1}2}\ls{3\sqrt 2}4\e \ls 9p$ $P_{[\f p4]}(-\f
79)$ $\mod p$. Thus applying Theorem 2.2 we have
$$\aligned P_{\f{p-1}2}\Ls{3\sqrt 2}4&\e \Ls{-6}p(\sqrt 2)^{\f{p-1}2}P_{[\f p4]}\Big(-\f
79\Big) \\&\e\cases 2^{\f{p-1}4}\cdot 2a\mod p&\t{if $p\e 1\mod 4$,
$p=a^2+4b^2$ and $4\mid a-1$,}
\\0\mod p&\t{if $p\e 3\mod 4$.}
\endcases\endaligned$$
It is well known that (see [BEW, Theorem 8.2.6])
$$2^{\f{p-1}4}\e \cases(-1)^{\f b2}\mod p&\t{if $p\e 1\mod 8$,}\\
\f {2b}a\mod p&\t{if $p\e 5\mod 8$, $p=a^2+4b^2$ and $b\e a\e 1\mod
4$.}
\endcases$$
Thus the result follows.

\pro{Theorem 2.7} Let $p$ be an odd prime. Then
$$P_{\f{p-1}2}(\sqrt 2)\e \cases
(-1)^{\f{p-1}8+\f{c-1}4}2c\mod p&\t{if $8\mid p-1$, $p=c^2+2d^2$ and
$c\e 1\mod 4$,}\\ -2\sqrt 2d\mod p&\t{if $8\mid p-3$, $p=c^2+2d^2$
and $d\e 1\mod 4$,}
\\0\mod p&\t{if $p\e 5,7\mod 8$.}
\endcases$$
\endpro
Proof. From [BE, Theorems 5.12 and 5.17] we know that
$$\sum_{k=0}^{p-1}\Ls{x^3-4x^2+2x}p=\cases (-1)^{[\f p8]+1}2c&\t{if $p=c^2+2d^2\e
1,3\mod 8$ with $4\mid c-1$,}
\\0&\t{otherwise.}\endcases$$
As $27(x^3-4x^2+2x)=(3x-4)^3-30(3x-4)-56$, we see that
$$\align&\sum_{x=0}^{p-1}\Ls{x^3-4x^2+2x}p
\\&=\Ls 3p\sum_{x=0}^{p-1}\Ls{(3x-4)^3-30(3x-4)-56}p =\Ls
3p\sum_{x=0}^{p-1}\Ls{x^3-30x-56}p .
\endalign$$ Thus, from the above we deduce

$$\aligned \sum_{x=0}^{p-1}\Ls{x^3-30x-56}p
= \cases (-1)^{\f{p+7}8}\sls p32c&\t{if $p\e 1\mod 8$, $p=c^2+2d^2$
and $4\mid c-1$,}
\\(-1)^{\f{p-3}8}\sls p32c&\t{if $p\e 3\mod 8$,
$p=c^2+2d^2$ and $4\mid c-1$,}
\\0&\t{if $p\e 5,7\mod 8$.}
\endcases\endaligned\tag 2.6$$
Taking $t=2$ in Theorem 2.1(i) and applying (2.6) we get
$$\aligned&\Ls 2p(\sqrt 2)^{\f{p-1}2}P_{\f{p-1}2}( \sqrt 2)\\&\e
P_{[\f p4]}(0)\e \sum_{k=0}^{[p/4]}\f{\b{4k}{2k}\b{2k}k}{128^k} \e
-\Ls{-12}p\sum_{x=0}^{p-1}\Ls{2x^3-15x-14}p
\\&=-\ls{-3}p\sum_{x=0}^{p-1}\Ls{2(\f x2)^3-15(\f x2)-14}p
=-\Ls p3\sum_{x=0}^{p-1}\Ls{x^3-30x-56}p
\\&\e
\cases (-1)^{\f{p-1}8}2c\mod p&\t{if $8\mid p-1$, $p=c^2+2d^2$ and
$4\mid c-1$,}\\-(-1)^{\f{p-3}8}2c\mod p&\t{if $8\mid p-3$,
$p=c^2+2d^2$ and $4\mid c-1$,}
\\0\mod p&\t{if $p\e 5,7\mod 8$.}
\endcases\endaligned$$
If $p\e 1\mod 8$ and $p=c^2+2d^2$ with $c\e 1\mod 4$ and
$d=2^{\alpha}d_0(2\nmid d_0)$, then $(\sqrt
2)^{\f{p-1}2}=2^{\f{p-1}4}\e (c/d)^{\f{p-1}2}\e \sls cp\sls dp=\sls
pc\sls p{d_0}=\sls{c^2+2d^2}c \sls{c^2+2d^2}{d_0}=\sls
2c=(-1)^{\f{c-1}4}\mod p$. If $p\e 3\mod 8$ and $p=c^2+2d^2$ with
$c\e d\e 1\mod 4$, then $\ls cd^2\e -2\e 2^{\f{p+1}2}\mod p$ and so
$c/d \e \pm 2^{\f{p+1}4}\mod p.$ As $\ls{c/d}p=\ls cp\ls
dp=\ls{c^2+2d^2}c\ls{c^2+2d^2}d=(-1)^{\f{c^2-1}8}=(-1)^{\f{p-3-2(d^2-1)}8}
=(-1)^{\f{p-3}8}$, we see that $c/d \e -(-1)^{\f{p-3}8}$ $
2^{\f{p+1}4}\mod p.$ Hence $(-1)^{\f{p-3}8}2c( \sqrt
2)^{-\f{p-1}2}\e -2^{\f{p+1}4}d\cdot 2(\sqrt 2)^{1-\f{p+1}2}=-
2\sqrt 2d\mod p.$ Now combining all of the above we obtain the
result.
\par\q
\newline{\bf Remark 2.1.} In [Su1], using Pell sequence Zhi-Wei Sun made a conjecture
related to Theorem 2.7.
 \pro{Theorem 2.8} Let $p>3$ be a prime.
Then
$$P_{\f{p-1}2}(\sqrt {-3})\e \cases (-1)^{\f{p-1}4}2A\mod p
&\t{if $12\mid p-1$, $p=A^2+3B^2$ and
 $4\mid A-1$,}\\(-1)^{\f{p-3}4}2B\sqrt {-3}\mod p&\t{if $12\mid p-7$, $p=A^2+3B^2$ and
 $4\mid B-1$,}\\0\mod p&\t{if $p\e 2\mod 3$}
 \endcases$$
 and
 $$P_{\f{p-1}2}\Ls{\sqrt 3}2\e \cases 2A\mod p&\t{if $12\mid p-1$, $p=A^2+3B^2$ and
 $4\mid A-1$,}\\-2\sqrt 3B\mod p&\t{if $12\mid p-7$, $p=A^2+3B^2$ and
 $4\mid B-1$,}\\0\mod p&\t{if $p\e 2\mod 3$.}
 \endcases$$
 \endpro
Proof. From Theorems 2.1 and 2.3 we have
$$\aligned &(\sqrt {-3})^{-\f{p-1}2}P_{\f{p-1}2}(\sqrt {-3})
\\&\e P_{[\f p4]}\Big(-\f 53\Big)\e \cases 2A\sls A3\mod p&\t{if $3\mid p-1$ and
$p=A^2+3B^2$,}\\0\mod p&\t{if $p\e 2\mod 3$.}
 \endcases\endaligned$$
 Thus the result is true for $p\e 2\mod 3$.
Now suppose $p\e 1\mod 3$, $p=A^2+3B^2$ and $A\e 1\mod 4$. By [S2,
p.1317] we have
$$3^{[\f p4]}\e \cases \sls A3\mod p&\t{if $p\e 1\mod{12}$,}
\\\sls A3\f BA\mod p&\t{if $p\e 7\mod{12}$ and $B\e 1\mod 4$.}\endcases\tag 2.7$$
Thus,
$$P_{\f{p-1}2}(\sqrt{-3})
\e\cases (-3)^{\f{p-1}4}\cdot 2A\sls A3\e (-1)^{\f{p-1}4}2A\mod
p\q\t{if $p\e 1\mod{12}$,}
\\(-3)^{\f{p-3}4}\sqrt{-3}\cdot 2A\sls A3\e (-1)^{\f{p-3}4}2B\sqrt{-3}\mod
p\\\qq\qq\qq\qq\qq\qq\qq\qq\qq\t{if $12\mid p-7$ and $4\mid
B-1$.}\endcases$$ Taking $x=-3$ in Lemma 2.5 we obtain $(\sqrt
3)^{\f{p-1}2}P_{\f{p-1}2}\sls{\sqrt 3}2\e (\sqrt{-3})^{\f{p-1}2}$
$P_{\f{p-1}2}(\sqrt{-3})$ $\mod p$. Now combining all the above we
obtain the result.

 \subheading{3.
Congruences for $\sum_{k=0}^{\f{p-1}2}\f 1{m^k}\b{2k}k^3\mod {p^2}$}

\pro{Lemma 3.1} Let $m$ and $n$ be nonnegative integers. Then
$$\sum_{k=0}^m\b{2k}k\b{n+k}{2k}\b{2k}k\b{k}{m-k}
=\sum_{k=0}^m\b{2k}k\b{n+k}{2k}\b{2m-2k}{m-k}\b{n+m-k}{2m-2k}.$$
\endpro
\par Lemma 3.1 can be easily proved by using WZ method. For the WZ method one may consult
[PWZ]. Clearly the lemma is true for $m=0,1$. Using Mathematica we
find both sides satisfy the same recurrence relation
$$\align (m+1)&(m+2n+2)(m-2n)f(m)+(2m+3)(m^2-2n^2
+3m-2n+2)f(m+1)\\&+(m+2)^3f(m+2)=0.\endalign$$ Thus the lemma is
true. The proof certificate for the left hand is
$$R_1(m,k)=-\f{k^2(2k-m-1)(2k-m)(m+2)}{(m-k+1)(m+2-k)},$$ and
the proof certificate for the right hand is
$$R_2(m,k)=\f{k^2(m-n-k)(m+n-k+1)P(m,k)}{(m-k+1)^2(m-k+2)^2},$$
where
$$P(m,k)=3mn^2-2n^2k+m^2+3mn-mk+6n^2-2nk+4m+6n-2k+4.$$

 \pro{Definition 3.1} Note that $\b nk\b{n+k}k=\b{2k}k\b{n+k}{2k}$.
For any nonnegative integer $n$ we define
$$S_n(x)=\sum_{k=0}^n\b nk\b{n+k}k\b{2k}kx^k=\sum_{k=0}^n
\b{2k}k^2\b{n+k}{2k}x^k.$$
\endpro

 \pro{Lemma 3.2} Let $n$ be a nonnegative
integer. Then
$$S_n(x)=P_n(\sqrt{1+4x})^2\qtq{and}P_n(x)^2=S_n\Ls{x^2-1}4.$$
\endpro
Proof. From (2.1) we have
$$\align P_n(x)^2&=\Big(\sum_{k=0}^n\b{2k}k\b{n+k}{2k}\Ls{x-1}2^k\Big)
\Big(\sum_{r=0}^n\b{2r}r\b{n+r}{2r}\Ls{x-1}2^r\Big)
\\&=\sum_{m=0}^{2n}\Ls{x-1}2^m\sum_{k=0}^m\b{2k}k\b{n+k}{2k}
\b{2m-2k}{m-k}\b{n+m-k}{2m-2k}.
\endalign$$
On the other hand,
$$\align
S_n\Ls{x^2-1}4&=\sum_{k=0}^n\b{2k}k^2\b{n+k}{2k}\Ls{x-1}2^k
\Big(1+\f{x-1}2\Big)^k
\\&=\sum_{k=0}^n\b{2k}k^2\b{n+k}{2k}\Ls{x-1}2^k\sum_{r=0}^k\b
kr\Ls{x-1}2^r
\\&=\sum_{m=0}^{2n}\Ls{x-1}2^m\sum_{k=0}^m\b{2k}k^2\b{n+k}{2k}
\b k{m-k}.
\endalign$$
Hence, from the above and Lemma 3.1 we deduce
$$P_n(x)^2=S_n\Ls{x^2-1}4\qtq{and so} P_n(\sqrt{1+4x})^2=S_n(x).$$
This proves the lemma.

 \pro{Theorem
3.1} Let $p$ be an odd prime and $m\in\Bbb Z_p$ with $m\not\e 0\mod
p$. Then
$$\sum_{k=0}^{\f{p-1}2}\f{\b{2k}k^3}{m^k}\e
P_{\f{p-1}2}\Big(\sqrt{1-\f{64}m}\Big)^2\mod {p^2}.$$
\endpro
Proof. By Definition 3.1 and (1.1) we have
$$S_{\f{p-1}2}\Big(-\f{16}m\Big)=\sum_{k=0}^{\f{p-1}2}\b{2k}k^2\b{\f{p-1}2+k}
{2k}\Big(-\f{16}m\Big)^k\e
\sum_{k=0}^{\f{p-1}2}\f{\b{2k}k^3}{m^k}\mod {p^2}.$$ On the other
hand, by Lemma 3.2 we get
$$S_{\f{p-1}2}\Big(-\f{16}m\Big)
=P_{\f{p-1}2}\Big(\sqrt{1+4\big(-\f{16}m\big)}\Big)^2=P_{\f{p-1}2}
\Big(\sqrt{\f{m-64}m}\Big)^2.$$ Thus the result follows.

\pro{Theorem 3.2} Let $p\not=2,7$ be a prime. Then
$$\sum_{k=0}^{\f{p-1}2}\b{2k}k^3\e (-1)^{\f{p-1}2}\sum_{k=0}^{\f{p-1}2}\f{\b{2k}k^3}{4096^k}
\e \cases 0\mod{p^2}&\t{if $p\e 3,5,6\mod 7$,}
\\4C^2\mod p&\t{if $p=C^2+7D^2\e 1,2,4\mod 7$.}\endcases$$
\endpro

Proof. Taking $m=1,4096$ in Theorem 3.1 and then applying Theorem
2.5 we deduce the result.
\newline{\bf Remark 3.1}
 The
congruence modulo $p$ can be deduced from (1.2) and [O, Corollary
11].

 \pro{Theorem 3.3} Let $p$ be an odd prime. Then
$$\sum_{k=0}^{\f{p-1}2}\f{\b{2k}k^3}{(-8)^k}
\e\cases 0\mod{p^2}&\t{if $p\e 3\mod 4$,}
\\4a^2-2p\mod{p^2}&\t{if $p=a^2+4b^2\e 1\mod 4$.}
\endcases$$
\endpro
Proof. From Theorem 3.1 we have
$\sum_{k=0}^{\f{p-1}2}\b{2k}k^3/(-8)^k \e
P_{\f{p-1}2}(3)^2\mod{p^2}.$ According to [S3, Corollary 2.3 and
Theorem 2.9], we have $$P_{\f{p-1}2}(3)\e\cases 0\mod p&\t{if $p\e
3\mod 4$,}\\(-1)^{\f{p-1}4}(2a-\f p{2a})\mod{p^2}&\t{if $p=a^2+4b^2
\e 1\mod 4$.}\endcases$$ Now combining all the above we deduce the
result.

\pro{Theorem 3.4} Let $p$ be an odd prime. Then
$$\align&\sum_{k=0}^{\f{p-1}2}\f{\b{2k}k^3}{(-512)^k}
\e  0\mod{p^2}\qtq{for $p\e 3\mod 4$,}\tag i
\\&\sum_{k=0}^{\f{p-1}2}\f{\b{2k}k^3}{(-64)^k}
\e  0\mod{p^2}\qtq{for $p\e 5,7\mod 8$,}\tag ii
\\&\sum_{k=0}^{\f{p-1}2}\f{\b{2k}k^3}{16^k}
\e  \sum_{k=0}^{\f{p-1}2}\f{\b{2k}k^3}{256^k}\e 0\mod{p^2}\qtq{for
$p\e 2\mod 3$.}\tag iii
\endalign$$
\endpro
Proof. Putting $m=-512,-64,16,256$ in Theorem 3.1 and then applying
Theorems 2.6-2.8 we deduce the result.
 \newline{\bf Remark 3.2} Theorems 3.3 and 3.4 were conjectured by Zhi-Wei Sun in
 [Su1,Su3].

\subheading{4. A general congruence modulo $p^2$}
\pro{Lemma 4.1}
Let $m$ be a nonnegative integer. Then
$$\sum_{k=0}^m\b {2k}k^2\b{4k}{2k}\b k{m-k}(-64)^{m-k}
=\sum_{k=0}^m\b {2k}k\b{4k}{2k}\b{2(m-k)}{m-k}\b{4(m-k)}{2(m-k)}.$$
\endpro
\par We prove the lemma by using WZ method and Mathematica. Clearly
the result is true for $m=0,1$. Since both sides satisfy the same
recurrence relation $$\align 1024&(m+1)(2m+1)(2m+3)S(m)
    -8(2m+3)(8m^2+24m+19)S(m+1)\\&+(m+2)^3 S(m+2) = 0,\endalign$$
    we see that
    Lemma 4.1 is true.
 The  proof certificate for the left hand side is
$$ - \frac{4096 k^2(m+2)(m-2k)(m-2k+1)}{(m-k+1)(m-k+2)},$$
 and the proof certificate for the right hand side is
$$ \frac{16
k^2(4m-4k+1)(4m-4k+3)(16m^2-16mk+55m-26k+46)}{(m-k+1)^2(m-k+2)^2}.$$
\pro{Theorem 4.1} Let $p$ be an odd prime and let  $x$ be a
variable. Then
$$\sum_{k=0}^{p-1}\b{2k}k^2\b{4k}{2k}(x(1-64x))^k\e \Big(
\sum_{k=0}^{p-1}\b{2k}k\b{4k}{2k}x^k\Big)^2\mod {p^2}.$$
\endpro
Proof. It is clear that
$$\align &\sum_{k=0}^{p-1}\b{2k}k^2\b{4k}{2k}(x(1-64x))^k
\\&=\sum_{k=0}^{p-1}\b{2k}k^2\b{4k}{2k}x^k\sum_{r=0}^k\b kr(-64x)^r
\\&=\sum_{m=0}^{2(p-1)}x^m\sum_{k=0}^{min\{m,p-1\}}\b{2k}k^2\b{4k}{2k}\b
k{m-k}(-64)^{m-k}.\endalign$$
 Suppose $p\le m\le 2p-2$ and $0\le
k\le p-1$. If $k>\f p2$, then $p\mid \b{2k}k$ and so $p^2\mid
\b{2k}k^2$. If $k<\f p2$, then $m-k\ge p-k>k$ and so $\b k{m-k}=0$.
Thus, from the above and Lemma 4.1 we deduce
$$\align &\sum_{k=0}^{p-1}\b{2k}k^2\b{4k}{2k}(x(1-64x))^k
\\&\e \sum_{m=0}^{p-1}x^m\sum_{k=0}^m\b {2k}k^2\b{4k}{2k}
\b k{m-k}(-64)^{m-k}
\\&=\sum_{m=0}^{p-1}x^m\sum_{k=0}^m\b {2k}k\b{4k}{2k}
\b{2(m-k)}{m-k}\b{4(m-k)}{2(m-k)}
\\&=\sum_{k=0}^{p-1}\b{2k}k\b{4k}{2k}x^k\sum_{m=k}^{p-1}
\b{2(m-k)}{m-k}\b{4(m-k)}{2(m-k)}x^{m-k}
\\&=\sum_{k=0}^{p-1}\b{2k}k\b{4k}{2k}x^k
\sum_{r=0}^{p-1-k}\b{2r}r\b{4r}{2r}x^r
\\&=\Big(\sum_{k=0}^{p-1}\b{2k}k\b{4k}{2k}x^k\Big)^2
-\sum_{k=0}^{p-1}\b{2k}k\b{4k}{2k}x^k
\sum_{r=p-k}^{p-1}\b{2r}r\b{4r}{2r}x^r
 \mod{p^2}.\endalign$$
Now suppose $0\le k\le p-1$ and $p-k\le r\le p-1$. If $k\ge
\f{3p}4$, then $p^2\nmid (2k)!$, $p^3\mid (4k)!$ and so $
\b{2k}k\b{4k}{2k}=\f{(4k)!}{(2k)!k!^2}\e 0\mod{p^2}$. If $k < \f
p4$, then $r\ge p-k> \f{3p}4$ and so $
\b{2r}r\b{4r}{2r}=\f{(4r)!}{(2r)!r!^2}\e 0\mod{p^2}$. If $\f p4<k<\f
p2$, then $r\ge p-k>\f p2$, $p\mid \b{2r}r$ and $p\mid \b{4k}{2k}$.
If $\f p2<k<\f{3p}4$, then $r\ge p-k>\f p4$, $p\mid \b{2k}k$ and
$p\mid \b{2r}r\b{4r}{2r}$. Hence we always have
$\b{2k}k\b{4k}{2k}\b{2r}r\b{4r}{2r}\e 0\mod{p^2}$ and so
$$\sum_{k=0}^{p-1}\b{2k}k\b{4k}{2k}x^k
\sum_{r=p-k}^{p-1}\b{2r}r\b{4r}{2r}x^r
 \e 0\mod{p^2}.$$
Now combining all the above we obtain the result.

 \pro{Corollary 4.1} Let $p$ be an odd
prime and $m\in\Bbb Z_p$ with $m\not\e 0\mod p$. Then
$$\sum_{k=0}^{p-1}\f{\b{2k}k^2\b{4k}{2k}}{m^k}
\e\Big(\sum_{k=0}^{p-1}\b{2k}k\b{4k}{2k}
\Big(\f{1-\sqrt{1-256/m}}{128}\Big)^k\Big)^2 \mod{p^2}.$$
\endpro
Proof. Taking $x=\f{1-\sqrt{1-256/m}}{128}$ in Theorem 4.1 we deduce
the result.

\pro{Lemma 4.2} Let $p$ be a prime greater than $3$. Then
$$P_{[\f p4]}(\sqrt t)\e -\sum_{x=0}^{p-1}(x^3+4x^2+2(1-\sqrt
t)x)^{\f{p-1}2}\mod p.$$
\endpro
Proof. It is clear that
$$\align&\sum_{x=0}^{p-1}(x^3+4x^2+2(1-\sqrt
t)x)^{\f{p-1}2}
\\&\e\sum_{x=0}^{p-1}\Big(\big(x-\f 43\big)^3+4\big(x-\f 43\big)^2+2(1-\sqrt
t)\big(x-\f 43\big)\Big)^{\f{p-1}2}
\\&=\sum_{x=0}^{p-1}\Big(x^3-\f{2(3\sqrt t+5)}3x+\f{8(9\sqrt
t+7)}{27}\Big)^{\f{p-1}2}
\\&\e \sum_{x=0}^{p-1}\Big(\ls{2x}3^3-\f{2(3\sqrt t+5)}3\cdot \f{2x}3
+\f{8(9\sqrt t+7)}{27}\Big)^{\f{p-1}2}
\\&=\Ls 8{27}^{\f{p-1}2}\sum_{x=0}^{p-1}\big(x^3-\f 32(3\sqrt
t+5)x+9\sqrt t+7\big)^{\f{p-1}2}\mod p.
\endalign$$
By Theorem 2.1(ii),
$$P_{[\f p4]}(u)\e -\ls 6p
\sum_{x=0}^{p-1}\big(x^3-\f 32(3u+5)x+9u+7\big)^{\f{p-1}2}\mod p\tag
4.1$$ for $u=0,1,\ldots,p-1$. Since both sides are polynomials of
$u$ with degree at most $(p-1)/2$. By Lagrange's theorem, (4.1) is
also true when $u$ is a variable. Hence
$$P_{[\f p4]}(\sqrt t)\e -\ls 6p
\sum_{x=0}^{p-1}\big(x^3-\f 32(3\sqrt t+5)x+9\sqrt
t+7\big)^{\f{p-1}2}\mod p.$$ Now combining all the above with the
fact $\sls 8{27}^{\f{p-1}2}\e \ls 6p\mod p$ we deduce the result.

 \pro{Theorem 4.2} Let $p$ be an odd prime, $m\in\Bbb Z_p$, $m\not\e 0\mod
p$ and $t=\sqrt{1-256/m}$. Then
$$\align \sum_{k=0}^{p-1}\f{\b{2k}k^2\b{4k}{2k}}{m^k}\e P_{[\f
p4]}(t)^2\e\Big(\sum_{x=0}^{p-1}(x^3+4x^2
+2(1-t)x)^{\f{p-1}2}\Big)^2\mod p.\endalign$$ Moreover, if $P_{[\f
p4]}(t)\e 0\mod p$ or $\sum_{x=0}^{p-1}(x^3+4x^2
+2(1-t)x)^{\f{p-1}2}\e 0\mod p$, then
$$\sum_{k=0}^{p-1}\f{\b{2k}k^2\b{4k}{2k}}{m^k}\e 0\mod{p^2}.$$
\endpro
Proof. For $\f p2<k<p$, $\b{2k}k=\f{(2k)!}{k!^2}\e 0\mod p$. For $\f
p4<k<\f p2$, $\b{4k}{2k}=\f{(4k)!}{(2k)!^2}\e 0\mod p$. Thus, $p\mid
\b{2k}k\b{4k}{2k}$ for $\f p4<k<p$. Hence, by Lemmas 2.2 and 4.2,
$$P_{[\f p4]}(t)\e
\sum_{k=0}^{p-1}\b{2k}k\b{4k}{2k}\Ls{1-t}{128}^k\e
-\sum_{x=0}^{p-1}(x^3+4x^2 +2(1-t)x)^{\f{p-1}2}\mod p.$$
 This together with Corollary 4.1 gives the result.

\pro{Theorem 4.3} Let $p\e 1,3\mod 8$ be a prime and so $p=c^2+2d^2$
with $c,d\in\Bbb Z$ and $c\e 1\mod 4$. Then
$$\sum_{k=0}^{p-1}\f{\b{2k}k\b{4k}{2k}}{128^k}
\e (-1)^{[\f p8]+\f{p-1}2}\Big(2c-\f p{2c}\Big)\mod{p^2}.$$
\endpro
Proof. By  the proof of Theorem 2.7,
$$\sum_{k=0}^{p-1}\f{\b{2k}k\b{4k}{2k}}{128^k}\e
\sum_{k=0}^{[p/4]}\f{\b{2k}k\b{4k}{2k}}{128^k}
 \e (-1)^{[\f p8]+\f{p-1}2}2c\mod p.$$
  Set $\sum_{k=0}^{p-1}\f{\b{2k}k\b{4k}{2k}}{128^k}
=(-1)^{[\f p8]+\f{p-1}2}2c+qp$. Then
 $$\Big(\sum_{k=0}^{p-1}\f{\b{2k}k\b{4k}{2k}}{128^k}\Big)^2
 =((-1)^{[\f p8]+\f{p-1}2}2c+qp)^2\e
 4c^2+(-1)^{[\f p8]+\f{p-1}2}4cqp\mod{p^2}.$$ Taking $x=\f 1{128}$
  in Theorem 4.1 we get
 $$\sum_{k=0}^{p-1}\f{\b{2k}k^2\b{4k}{2k}}{256^k}
 \e \Big(\sum_{k=0}^{p-1}\f{\b{2k}k\b{4k}{2k}}{128^k}\Big)^2\mod{p^2}.$$
 From [M] and [Su2] we have
 $\sum_{k=0}^{p-1}\b{2k}k^2\b{4k}{2k}/256^k\e 4c^2-2p\mod{p^2}.$
 Thus
 $$4c^2-2p\e  \Big(\sum_{k=0}^{p-1}\f{\b{2k}k\b{4k}{2k}}{128^k}\Big)^2
 \e 4c^2+(-1)^{[\f p8]+\f{p-1}2}4cqp\mod{p^2}$$ and hence
 $q\e -(-1)^{[\f p8]+\f{p-1}2}\f 1{2c}\mod p$. So the theorem is proved.
\par\q
\par We note that Theorem 4.3 was conjectured by Zhi-Wei Sun in
[Su1].

\subheading{5. Congruences for $\sum_{k=0}^{p-1}
 \b{2k}k^2\b{4k}{2k}\big/m^k$}

\pro{Theorem 5.1} Let $p\not=2,3,7$ be a prime. Then
$$\aligned &\sum_{k=0}^{p-1}\f{\b{2k}k^2\b{4k}{2k}}{648^k}\e
\cases 0\mod {p^2}&\t{if $p\e 3\mod 4$,}
\\4a^2\mod p&\t{if $p=a^2+4b^2\e 1\mod 4$,}
\endcases
\\&\sum_{k=0}^{p-1}\f{\b{2k}k^2\b{4k}{2k}}{(-144)^k}\e
\cases 0\mod {p^2}&\t{if $p\e 2\mod 3$,}
\\4A^2\mod p&\t{if $p=A^2+3B^2\e 1\mod 3$,}
\endcases
\\&\sum_{k=0}^{p-1}\f{\b{2k}k^2\b{4k}{2k}}{(-3969)^k}\e
\cases 0\mod {p^2}&\t{if $p\e 3,5,6\mod 7$,}
\\4C^2\mod p&\t{if $p=C^2+7D^2\e 1,2,4\mod 7$.}
\endcases\endaligned$$
\endpro
Proof. Taking $m=648,-144,-3969$ in Theorem 4.2 and then applying
Theorems 2.2-2.4 and (1.6) we deduce the result.
\par We remark that Theorem 5.1 was conjectured by the author in
[S3].
 \pro{Lemma 5.1 ([S4, Lemma 4.1])} Let $p$ be an odd prime and
let $a,m,n$ be p-adic integers. Then
$$\sum_{x=0}^{p-1}(x^3+a^2mx+a^3n)^{\f{p-1}2}\e a^{\f{p-1}2}
\sum_{x=0}^{p-1}(x^3+mx+n)^{\f{p-1}2}\mod p.$$ Moreover, if $a,m,n$
are congruent to some integers modulo $p$, then
$$\sum_{x=0}^{p-1}\Ls{x^3+a^2mx+a^3n}p=
\Ls ap\sum_{x=0}^{p-1}\Ls{x^3+mx+n}p.$$
\endpro

\pro{Theorem 5.2} Let $p\not=2,3,7$ be a prime. Then
$$\aligned &P_{[\f p4]}\Ls{5\sqrt{-7}}9\\&\e\cases
\sls {3(7+\sqrt{-7})}p\sls C72C\mod p&\t{if $p\e 1,2,4\mod 7$ and so
$p=C^2+7D^2$,}
\\0\mod p&\t{if $p\e 3,5,6\mod 7$}
\endcases\endaligned$$
and
$$\sum_{k=0}^{p-1}\f{\b{2k}k^2\b{4k}{2k}}{81^k}
\e\cases 4C^2\mod p &\t{if $p\e 1,2,4\mod 7$ and so $p=C^2+7D^2$,}
\\0\mod {p^2}&\t{if $p\e 3,5,6\mod 7$.}
\endcases$$
\endpro
Proof. By (4.1),
$$P_{[\f p4]}\Ls{5\sqrt{-7}}9
\e -\Ls 6p\sum_{x=0}^{p-1}\Big(x^3-\f
52(3+\sqrt{-7})x+7+5\sqrt{-7}\Big)^{\f{p-1}2}\mod p.$$ Since
$$\f{-\f
52(3+\sqrt{-7})}{-35}=\Ls{1-\sqrt{-7}}{2\sqrt{-7}}^2\qtq{and}
\f{7+5\sqrt{-7}}{-98}=\Ls{1-\sqrt{-7}}{2\sqrt{-7}}^3,$$ by the above
and Lemma 5.1 we have
$$P_{[\f p4]}\Ls{5\sqrt{-7}}9
\e -\Ls
6p\Ls{1-\sqrt{-7}}{2\sqrt{-7}}^{\f{p-1}2}\sum_{x=0}^{p-1}\Ls{x^3-35x-98}p
\mod p.$$ By the proof of Theorem 2.4,
$$\aligned \sum_{x=0}^{p-1}\Ls{x^3-35x-98}p=\cases -2C\sls C7&\t{if $p=C^2+7D^2\e 1,2,4\mod 7$,}
\\0&\t{if $p\e 3,5,6\mod 7$.}
\endcases\endaligned\tag 5.1$$
For $p\e 1,2,4\mod 7$ we see that
$$\Ls 6p\Ls{1-\sqrt{-7}}{2\sqrt{-7}}^{\f{p-1}2}
=\Ls 6p\Ls{7+\sqrt{-7}}{2\cdot (-7)}^{\f{p-1}2}\e \Ls
3p\Ls{7+\sqrt{-7}}p\mod p.$$ Thus, from the above we deduce the
congruence for $P_{[\f p4]}\sls{5\sqrt{-7}}9\mod p$. Applying
Theorem 4.2 (with $m=81$) we obtain the remaining result.

\pro{Theorem 5.3} Let $p$ be a prime such that $p\e 1,9\mod{20}$ and
hence $p=u^2+5v^2$ for some integers $u$ and $v$. Then
$$\sum_{k=0}^{p-1}\f{\b{2k}k^2\b{4k}{2k}}{(-1024)^k}
\e 4u^2\mod p.$$
\endpro
Proof. Taking $m=-1024$ and $t=\sqrt 5/2$ in Theorem 4.2 we see that
$$\sum_{k=0}^{p-1}\f{\b{2k}k^2\b{4k}{2k}}{(-1024)^k}\e
\Big(\sum_{x=0}^{p-1}\Ls{x^3+4x^2+(2-\sqrt 5)x}p\Big)^2\mod p.$$ By
[LM, Theorem 11] we have
$$\sum_{x=0}^{p-1}\Ls{x^3+4x^2+(2-\sqrt
5)x}p=\pm 2u.$$ Thus the result follows.
\newline{\bf Remark 5.1}  Let $p\e 1\mod{20}$ be a prime and hence
$p=a^2+4b^2=u^2+5v^2$ with $a,b,u,v\in\Bbb Z$. A result of Cauchy
([BEW, p. 291]) states that
$$\b{\f{p-1}2}{\f{p-1}{20}}^2
\e\cases 4u^2\mod p&\t{if $5\nmid a$,}
\\-4u^2\mod p&\t{if $5\mid a$.}
\endcases$$ Let $m\in\{5,13,37\}$ and
$f(m)=-1024,-82944,-2^{10}\cdot 21^4$ according as $m=5,13$ or $37$.
Suppose that $p$ is an odd prime such that $p\nmid mf(m)$. In [Su1],
Z.W. Sun conjectured that
$$\sum_{k=0}^{p-1}\f{\b{2k}k^2\b{4k}{2k}}{f(m)^k}
\e\cases 4x^2-2p\mod {p^2}&\t{if $\sls {m}p=\sls{-1}p=1$ and so
$p=x^2+my^2$,}\\2p-2x^2\mod {p^2}&\t{if $\sls mp=\sls{-1}p=-1$ and
so $2p=x^2+my^2$,}
\\0\mod{p^2}&\t{if $\sls mp=-\sls{-1}p$.}
\endcases$$

\par\q
\par Let $p>3$ be a prime and let $\Bbb F_p$ be the field of $p$ elements.
For $m,n\in \Bbb F_p$ let $\#E_p(x^3+mx+n)$ be the number of points
on the curve $E_p$: $y^2=x^3+mx+n$ over the field $\Bbb F_p$. It is
well known that
$$\#E_p(x^3+mx+n)=p+1+\sum_{x=0}^{p-1}\Ls{x^3+mx+n}p.\tag 5.2$$
Let $K=\Bbb Q(\sqrt{-d})$ be an imaginary quadratic field and the
curve $y^2=x^3+mx+n$ has complex multiplication by $K$. By Deuring's
theorem ([C, Theorem 14.16],[PV],[I]), we have
$$\#E_p(x^3+mx+n)=\cases p+1&\t{if $p$ is inert in $K$,}
\\p+1-\pi-\bar{\pi}&\t{if $p=\pi\bar{\pi}$ in $K$,}
\endcases\tag 5.3$$ where $\pi$ is in an order in $K$ and
$\bar{\pi}$ is the conjugate number of $\pi$. If $4p=u^2+dv^2$ with
$u,v\in\Bbb Z$, we may take $\pi=\f 12(u+v\sqrt{-d})$. Thus,
$$\sum_{x=0}^{p-1}\Ls{x^3+mx+n}p=\cases \pm u&\t{if $4p=u^2+dv^2$
with $u,v\in\Bbb Z$,}\\0&\t{otherwise.}
\endcases\tag 5.4$$
In [JM] and [PV] the sign of $u$ in (5.4) was determined for those
imaginary quadratic fields $K$ with class number $1$. In [LM] and
[I] the sign of $u$ in (5.4) was determined for imaginary quadratic
fields $K$ with class number $2$.

\pro{Theorem 5.4} Let $p$ be a prime such that $p\e \pm 1\mod{12}$.
Then
$$P_{[\f p4]}\Big(\f 7{12}\sqrt 3\Big)\e
\cases\sls{2+2\sqrt 3}p2x
\mod p&\t{if $p=x^2+9y^2\e 1\mod{12}$ with $3\mid x-1$,} \\
0\mod p&\t{if $p\e 11\mod{12}$}\endcases
$$ and
$$\sum_{k=0}^{p-1}\f{\b{2k}k^2\b{4k}{2k}}{(-12288)^k}
\e\cases 4x^2\mod p&\t{if $p=x^2+9y^2\e 1\mod{12}$,}
\\0\mod{p^2}&\t{if $p\e 11\mod{12}$.}
\endcases$$
\endpro
Proof. From [I, p.133] we know that the elliptic curve defined by
the equation $y^2=x^3-(120+42\sqrt 3)x+448+336\sqrt 3$ has complex
multiplication by the order of discriminant $-36$. Thus, by (5.4)
and [I, Theorem 3.1] we have
$$\aligned&\sum_{n=0}^{p-1}\Ls{n^3-(120+42\sqrt 3)x+448+336\sqrt 3}
p\\&=\cases -2x\sls{1+\sqrt 3}p&\t{if $p=x^2+9y^2\e 1\mod{12}$ with
$3\mid x-1$,}
\\0&\t{if $p\e 11\mod{12}.$}\endcases\endaligned$$
By (4.1),
$$\align P_{[\f p4]}\Big(\f 7{12}\sqrt 3\Big)&\e
-\Ls 6p\sum_{n=0}^{p-1}\Big(n^3-\f{60+21\sqrt 3)}8n+\f{28+21\sqrt
3}4\Big)^{\f{p-1}2}
\\&\e -\Ls 6p\sum_{n=0}^{p-1}\Big(\ls n4^3-\f{60+21\sqrt 3)}8\cdot\f
n4+\f{28+21\sqrt 3}4\Big)^{\f{p-1}2}
\\&\e -\Ls 6p\sum_{n=0}^{p-1}\Ls{n^3-(120+42\sqrt 3)x+448+336\sqrt 3}
p\mod p.\endalign$$ Now combining all the above we obtain the
congruence for $P_{[\f p4]}(\f 7{12}\sqrt 3)\mod p$. Applying
Theorem 4.2 we deduce the remaining result.
\par\q
\newline{\bf Remark 5.2} In [Su1, Conjecture A24], Z.W. Sun
conjectured that for any prime $p>3$,
$$\align&\sum_{k=0}^{p-1}\f{\b{2k}k^2\b{4k}{2k}}{(-12288)^k}
\\&\e\cases (-1)^{[\f x6]}(4x^2-2p)\mod p&\t{if $p=x^2+y^2\e 1\mod{12}$
and $4\mid x-1$,}
\\-4\sls{xy}3xy\mod{p^2}&\t{if $p=x^2+y^2\e 5\mod{12}$ and $4\mid
x-1$,}
\\0\mod{p^2}&\t{if $p\e 3\mod{4}$.}
\endcases\endalign$$

 \pro{Theorem 5.5}
Let $p$ be an odd prime such that $p\not=3$ and $\sls {13}p=1$. Then
$$\sum_{k=0}^{p-1}\f{\b{2k}k^2\b{4k}{2k}}{(-82944)^k}
\e\cases 4x^2\mod p&\t{if $p=x^2+13y^2\e 1\mod 4$,}
\\0\mod{p^2}&\t{if $p\e 3\mod 4$.}
\endcases$$
\endpro
Proof. From [LM, Table II] we know that the elliptic curve defined
by the equation $y^2=x^3+4x^2+(2-\f 59\sqrt{13})x$ has complex
multiplication by the order of discriminant $-52$. Thus, by (5.4) we
have
$$\aligned\sum_{n=0}^{p-1}\Ls{n^3+4n^2+(2-\f 59\sqrt{13})n}p
=\cases 2x&\t{if $p\e 1\mod 4$ and so $p=x^2+13y^2$,}
\\0&\t{if $p\e 3\mod 4.$}\endcases\endaligned$$
Now taking $m=-2^{10}\cdot 3^4$ and $t=\f 5{18}\sqrt{13}$ in Theorem
4.2 and applying the above we deduce the result.

\pro{Theorem 5.6} Let $p$ be an odd prime such that $p\not=3,7$ and
$\sls {37}p=1$. Then
$$\sum_{k=0}^{p-1}\f{\b{2k}k^2\b{4k}{2k}}{(-2^{10}\cdot 21^4)^k}
\e\cases 4x^2\mod p&\t{if $p\e 1\mod 4$ and so $p=x^2+37y^2$,}
\\0\mod{p^2}&\t{if $p\e 3\mod 4$.}
\endcases$$
\endpro
Proof. From [LM, Table II] we know that the elliptic curve defined
by the equation $y^2=x^3+4x^2+(2-\f {145}{441}\sqrt{37})x$ has
complex multiplication by the order of discriminant $-148$. Thus, by
(5.4) we have
$$\aligned\sum_{n=0}^{p-1}\Ls{n^3+4n^2+(2-\f {145}{441}\sqrt{37})n}p
=\cases 2x&\t{if $p\e 1\mod 4$ and so $p=x^2+37y^2$,}
\\0&\t{if $p\e 3\mod 4.$}\endcases\endaligned$$
Now taking $m=-2^{10}\cdot 21^4$ and $t=\f {145}{882}\sqrt{37}$ in
Theorem 4.2 and applying the above we deduce the result.
\par\q
\par Let $b\in\{3,5,11,29\}$ and $f(b)=48^2,12^4,1584^2,396^4$ according as
$b=3,5,11,29$. For any odd prime $p$ with $p\nmid bf(b)$, Z.W. Sun
conjectured that ([Su1, Conjectures A14, A16, A18 and A21])
$$\sum_{k=0}^{p-1}\f{\b{2k}k^2\b{4k}{2k}}{f(b)^k}
\e\cases 4x^2-2p\mod {p^2}&\t{if $\sls 2p=\sls{-b}p=1$ and so
$p=x^2+2by^2$,}\\2p-8x^2\mod {p^2}&\t{if $\sls 2p=\sls{-b}p=-1$ and
so $p=2x^2+by^2$,}
\\0\mod{p^2}&\t{if $\sls 2p=-\sls{-b}p$.}
\endcases\tag 5.5$$
\par Now we partially solve the above conjecture.
\pro{Theorem 5.7} Let $p$ be an odd prime such that $p\e \pm 1 \mod
8$. Then
$$P_{[\f p4]}\Ls{2\sqrt 2}3\e
\cases (-1)^{\f{p-1}2}\sls{\sqrt 2}p\sls x32x\mod p&\t{if
$p=x^2+6y^2\e 1,7\mod {24}$,}
\\0\mod p&\t{if $p\e 17,23\mod {24}$}\endcases $$ and
$$\sum_{k=0}^{p-1}\f{\b{2k}k^2\b{4k}{2k}}{48^{2k}}
\e\cases 4x^2\mod p&\t{if $p=x^2+6y^2\e 1,7\mod {24}$,}
\\0\mod{p^2}&\t{if $p\e 17,23\mod {24}$.}
\endcases$$
\endpro
Proof. From [I, p.133] we know that the elliptic curve defined by
the equation $y^2=x^3+(-21+12\sqrt 2)x-28+22\sqrt 2$ has complex
multiplication by the order of discriminant $-24$. Thus, by (5.4)
and [I, Theorem 3.1] we have
$$\aligned&\sum_{n=0}^{p-1}\Ls{n^3+(-21+12\sqrt 2)n-28+22\sqrt 2}p
\\&=\cases 2x\sls{2x}3\sls{1+\sqrt 2}p
&\t{if $p\e 1,7\mod {24}$ and so $p=x^2+6y^2$,}
\\0&\t{if $p\e 17,23\mod {24}.$}\endcases\endaligned$$
By (4.1),
$$\align P_{[\f p4]}\Ls{2\sqrt 2}3&\e
-\Ls 6p\sum_{n=0}^{p-1}\Big(n^3-\f {15+6\sqrt 2}2n+7+6\sqrt
2\Big)^{\f{p-1}2}\mod p.\endalign$$ Since
$$\f{-(15+6\sqrt 2)/2}
{-21+12\sqrt 2}=\Ls{\sqrt 2+1}{\sqrt 2}^2\qtq{and} \f{7+6\sqrt
2}{-28+22\sqrt 2}=\Ls{\sqrt 2+1}{\sqrt 2}^3,$$ by Lemma 5.1 and the
above we have
$$\align P_{[\f p4]}\Ls{2\sqrt 2}3&\e
-\Ls 6p  \Ls{\sqrt 2(\sqrt
2+1)}p\sum_{n=0}^{p-1}\Ls{n^3+(-21+12\sqrt 2)n-28+22\sqrt 2}p
\\&\e\cases -\sls 6p\sls{\sqrt 2}p2x\sls{2x}3
\mod p&\t{if $p=x^2+6y^2\e 1,7\mod{24}$,}
\\0\mod p&\t{if $p\e 17,23\mod p$.}
\endcases\endalign$$
This yields the result for $P_{[\f p4]}\sls{2\sqrt 2}3\mod p$.
Taking $m=48^2$ and $t=\f 23\sqrt 2$ in Theorem 4.2 and applying the
above we deduce the remaining result.

\pro{Theorem 5.8} Let $p$ be a prime such that $p\e \pm 1\mod 5$.
Then
$$\sum_{k=0}^{p-1}\f{\b{2k}k^2\b{4k}{2k}}{12^{4k}}
\e\cases 4x^2\mod p&\t{if $p=x^2+10y^2\e 1,9,11,19\mod {40}$,}
\\0\mod{p^2}&\t{if $p\e 21,29,31,39\mod {40}$.}
\endcases$$
\endpro
Proof. From [LM, Table II] we know that the elliptic curve defined
by the equation $y^2=x^3+4x^2+(2-\f 89\sqrt{5})x$ has complex
multiplication by the order of discriminant $-40$. Thus, by (5.4) we
have
$$\aligned\sum_{n=0}^{p-1}\Ls{n^3+4n^2+(2-\f 89\sqrt 5)n}p
=\cases 2x&\t{if $p\e 1,9,11,19\mod {40}$ and so $p=x^2+10y^2$,}
\\0&\t{if $p\e 21,29,31,39\mod {40}.$}\endcases\endaligned$$
Now taking $m=12^4$ and $t=\f 49\sqrt 5$ in Theorem 4.2 and applying
the above we deduce the result.
\par\q

\pro{Theorem 5.9} Let $p$ be a prime such that $p\e \pm 1\mod 8$.
Then
$$\sum_{k=0}^{p-1}\f{\b{2k}k^2\b{4k}{2k}}{1584^{2k}}
\e\cases 4x^2\mod p&\t{if $\sls p{11}=1$ and so $p=x^2+22y^2$,}
\\0\mod{p^2}&\t{if $\sls p{11}=-1$.}
\endcases$$
\endpro
Proof. From [LM, Table II] we know that the elliptic curve defined
by the equation $y^2=x^3+4x^2+(2-\f {140}{99}\sqrt{2})x$ has complex
multiplication by the order of discriminant $-88$. Thus, by (5.4) we
have
$$\aligned\sum_{n=0}^{p-1}\Ls{n^3+4n^2+(2-\f {140}{99}\sqrt{2})n}p
=\cases 2x&\t{if $\sls p{11}=1$ and so $p=x^2+22y^2$,}
\\0&\t{if $\sls p{11}=-1.$}\endcases\endaligned$$
Now taking $m=1584^2$ and $t=\f {70}{99}\sqrt 2$ in Theorem 4.2 and
applying the above we deduce the result.
 \pro{Theorem 5.10}
Let $p$ be an odd prime such that $\sls {29}p=1$. Then
$$\sum_{k=0}^{p-1}\f{\b{2k}k^2\b{4k}{2k}}{396^{4k}}
\e\cases 4x^2\mod p&\t{if $p=x^2+58y^2\e 1,3\mod 8$,}
\\0\mod{p^2}&\t{if $p\e 5,7\mod 8$.}
\endcases$$
\endpro
Proof. From [LM, Table II] we know that the elliptic curve defined
by the equation $y^2=x^3+4x^2+(2-\f {3640}{9801}\sqrt{29})x$ has
complex multiplication by the order of discriminant $-232$. Thus, by
(5.4) we have
$$\aligned\sum_{n=0}^{p-1}\Ls{n^3+4n^2+(2-\f {3640}{9801}\sqrt{29})n}p
=\cases 2x&\t{if $p\e 1,3\mod 8$ and so $p=x^2+58y^2$,}
\\0&\t{if $p\e 5,7\mod 8.$}\endcases\endaligned$$
Now taking $m=396^4$ and $t=\f {1820}{9801}\sqrt {29}$ in Theorem
4.2  and applying the above we deduce the result.

\pro{Theorem 5.11} Let $p$ be an odd prime such that $p\e
1,5,19,23\mod{24}$. Then
$$\sum_{k=0}^{p-1}\f{\b{2k}k^2\b{4k}{2k}}{28^{4k}}
\e\cases 4x^2\mod p&\t{if $p\e 1,19\mod{24}$ and so $p=x^2+18y^2$,}
\\0\mod{p^2}&\t{if $p\e 5,23\mod {24}$.}
\endcases$$
\endpro
Proof. From [LM, Table II] we know that the elliptic curve defined
by the equation $y^2=x^3+4x^2+(2-\f {40}{49}\sqrt{6})x$ has complex
multiplication by the order of discriminant $-72$. Thus, by (5.4)
 we have
$$\aligned\sum_{n=0}^{p-1}\Ls{n^3+4n^2+(2-\f {40}{49}\sqrt{6})n}p
=\cases 2x&\t{if $p\e 1,19\mod {24}$ and so $p=x^2+18y^2$,}
\\0&\t{if $p\e 5,23\mod {24}.$}\endcases\endaligned$$
Now taking $m=28^4$ and $t=\f {20}{49}\sqrt {6}$ in Theorem 4.2 and
applying the above we deduce the result.
\par\q
\newline{\bf Remark 5.3} Let $p\not=2,7$ be a prime. Z.W. Sun conjectured
that ([Su1, Conjecture A28])
$$\sum_{k=0}^{p-1}\f{\b{2k}k^2\b{4k}{2k}}{28^{4k}}
\e\cases 4x^2-2p\mod p&\t{if $p=x^2+2y^2\e 1,3\mod 8$,}
\\0\mod{p^2}&\t{if $p\e 5,7\mod 8$.}
\endcases$$

\pro{Theorem 5.12} Let $p$ be an odd prime such that $p\e \pm 1\mod
5$. Then
$$\sum_{k=0}^{p-1}\f{\b{2k}k^2\b{4k}{2k}}{(-2^{14}\cdot 3^4\cdot 5)^k}
\e\cases 4x^2\mod p&\t{if $p=x^2+25y^2$,}
\\0\mod{p^2}&\t{if $p\e 3\mod 4$.}
\endcases$$
\endpro
Proof. From [LM, Table II] we know that the elliptic curve defined
by the equation $y^2=x^3+4x^2+(2-\f {161}{180}\sqrt{5})x$ has
complex multiplication by the order of discriminant $-100$. Thus, by
(5.4)  we have
$$\aligned\sum_{n=0}^{p-1}\Ls{n^3+4n^2+(2-\f {161}{180}\sqrt{5})n}p
=\cases 2x&\t{if $p=x^2+25y^2$,}
\\0&\t{if $p\e 3\mod 4.$}\endcases\endaligned$$
Now taking $m=-2^{14}\cdot 3^4\cdot 5$ and $t=\f {161}{360}\sqrt{5}$
in Theorem 4.2 and applying the above we deduce the result.
\par\q
\newline{\bf Remark 5.4} Let $p>5$ be a prime. Z.W. Sun made a conjecture
([Su1, Conjecture A25]) equivalent to
$$\sum_{k=0}^{p-1}\f{\b{2k}k^2\b{4k}{2k}}{(-2^{14}\cdot 3^4\cdot 5)^k}
\e\cases 4x^2-2p\mod p&\t{if $p=x^2+25y^2$,}
\\-4xy\mod{p^2}&\t{if $p=x^2+y^2$ with $5\mid x-y$,}
\\0\mod{p^2}&\t{if $p\e 3\mod 4$.}
\endcases$$

\par\q
\newline{\bf Acknowledgements} The author is indebted to Prof. Qing-Hu
Hou at Nankai University for his help in finding recurrence
relations and proof certificates concerning Lemmas 3.1 and 4.1.

\enddocument

\Refs \widestnumber\key {BEW}
  \ref\key A\by S. Ahlgren\paper
  Gaussian hypergeometric series and combinatorial congruences
  \jour in:
Symbolic computation, number theory, special functions, physics and
combina- torics (Gainesville, FI, 1999), pp. 1-12, Dev. Math., Vol.
4, Kluwer, Dordrecht, 2001\endref

\ref\key BE\by B. C. Berndt and R. J. Evans\paper  Sums of Gauss
Eisenstein, Jacobi, Jacobsthal and Brewer\jour Illinois J. Math.\vol
23\yr 1979\pages 374-437\endref

 \ref \key BEW\by  B.C. Berndt, R.J. Evans and K.S.
Williams\book  Gauss and Jacobi Sums\publ John Wiley $\&$
Sons\publaddr New York\yr 1998\endref

\ref\key B\by F. Beukers \paper Another congruence for the Ap\'ery
numbers\jour J. Number Theory \vol 25 \yr 1987\pages 201-210\endref

\ref\key C\by D.A. Cox \book Primes of the Form $x^2+ny^2$: Fermat,
Class Field Theory, and Complex Multiplication\publ Wiley\publaddr
New York\yr 1989\endref

 \ref \key G\by H.W. Gould\book Combinatorial
Identities, A Standardized Set of Tables Listing 500 Binomial
Coefficient Summations\publ Morgantown, W. Va.\yr 1972\endref

\ref\key I\by N. Ishii\paper Trace of Frobenius endomorphism of an
elliptic curve with complex multiplication \jour Bull. Austral.
Math. Soc.\vol 70\yr 2004 \pages 125-142\endref

\ref\key Is\by  T. Ishikawa\paper Super congruence for the Ap¡äery
numbers\jour Nagoya Math. J. \vol 118\yr 1990\pages 195-202\endref

 \ref\key JM\by A. Joux et F. Morain\paper Sur
les sommes de caract$\grave e$res li\'ees aux courbes elliptiques
$\grave a$ multiplication complexe \jour J. Number Theory\vol 55\yr
1995\pages 108-128\endref

 \ref\key LM\by F.
Lepr$\acute {\t{e}}$vost and F. Morain \paper Rev$\Hat
{\t{e}}$tements de courbes elliptiques $\grave {\t{a}}$
multiplication complexe par des courbes hyperelliptiques et sommes
de caract$\grave {\t{e}}$res \jour J. Number Theory \vol 64\yr
1997\pages 165-182\endref

\ref\key LR\by P.-R. Loh and R.C. Rhoades \paper p-adic and
combinatorial properties of modular form coefficients\jour Int. J.
Number Theory \vol 2\yr 2006\pages 305-328\endref

 \ref\key MOS\by W. Magnus, F. Oberhettinger and R.P. Soni\book Formulas and Theorems
for the Special Functions of Mathematical Physics, 3rd. ed.\publ
Springer\publaddr New York\yr 1966\pages 228-232\endref

 \ref\key M\by  E. Mortenson\paper Supercongruences for truncated $\
_{n+1}F_n$ hypergeometric series with applications to certain weight
three newforms\jour Proc. Amer. Math. Soc.\vol 133\yr 2005\pages
321-330.\endref

  \ref\key O\by K. Ono\paper Values of Gaussian hypergeometric
series\jour Trans. Amer. Math. Soc. \vol 350\yr 1998\pages
1205-1223\endref

 \ref\key
PV\by R. Padma and S. Venkataraman\paper Elliptic curves with
complex multiplication and a character sum\jour J. Number Theory\vol
61\yr 1996\pages 274-282\endref

\ref\key PWZ\by M. Petkov$\check {\t{s}}$ek, H. S. Wilf and D.
Zeilberger\book  A = B\publ A K Peters, Wellesley\endref

 \ref\key
R1\by A.R. Rajwade \paper The Diophantine equation
$y^2=x(x^2+21Dx+112D^2)$ and the conjectures of Birch and
Swinnerton-Dyer \jour J. Austral. Math. Soc. Ser. A\vol 24\yr 1977
\pages 286-295\endref

\ref\key R2\by A.R. Rajwade \paper On a conjecture of Williams\jour
Bull. Soc. Math. Belg. Ser. B\vol 36\yr 1984\pages 1-4\endref

 \ref\key
S1\by Z.H. Sun\paper On the number of incongruent residues of
$x^4+ax^2+bx$ modulo $p$\jour J. Number Theory \vol 119\yr
2006\pages 210-241\endref

 \ref\key S2\by Z.H. Sun\paper On the
quadratic character of quadratic units\jour J. Number Theory \vol
128\yr 2008\pages 1295-1335\endref


\ref\key S3\by Z.H. Sun\paper Congruences concerning Legendre
polynomials\jour Proc. Amer. Math. Soc.\vol 139\yr 2011\pages
1915-1929\endref

\ref\key S4\by Z.H. Sun\paper Congruences involving
$\b{2k}k^2\b{3k}km^{-k}$, arXiv:1104.2789.
http://arxiv.org/abs/1104.2789\endref


 \ref \key Su1\by Z.W. Sun\paper
Open conjectures on congruences, arXiv:0911.5665.
http://arxiv.org/abs/0911.5665\endref

\ref \key Su2\by Z.W. Sun\paper On sums involving products of three
binomial coefficients, preprint, arXiv:1012.3141.
http://arxiv.org/abs/1012.3141\endref


\ref \key Su3\by Z.W. Sun\paper On congruences related to central
binomial coefficients\jour J. Number Theory \vol 131\yr 2011\pages
2219-2238\endref

\ref \key Su4\by Z.W. Sun\paper Super congruences and Euler
numbers\jour Sci. China Math.\vol 54\yr 2011\pages 2509-2535\endref


\ref\key vH \by L. van Hamme\paper Some conjectures concerning
partial sums of generalized hypergeometric series \jour in: p-adic
Functional Analysis (Nijmegen, 1996), pp. 223-236, Lecture Notes in
Pure and Appl. Math., Vol. 192, Dekker, 1997\endref
\endRefs
\bye